\documentclass[11pt,amsfonts]{article}
\usepackage{graphicx}
\usepackage{latexsym}
\usepackage{amssymb}
\usepackage{amsmath}
\usepackage{enumerate}
\usepackage{color}
\usepackage{layout}
\usepackage{dsfont}
\usepackage{eufrak}

\newtheorem{lemma}{Lemma}

\newtheorem{corollary}{Corollary}
\newtheorem{theorem}{Theorem}
\newtheorem{remark}{Remark}

\def\real{{\mathord{{\rm I\kern-2.8pt R}}}}        % Fake blackboard bold R.
\def\inte{{\mathord{{\rm I\kern-2.8pt N}}}}

\def\sZZ{{\rm Z\kern-2.8ptem{}Z}}

\def\z{{\mathchoice
  {\sZZ}
  {\sZZ}
  {\rm Z\kern-0.30em{}Z}
  {\rm Z\kern-0.25em{}Z} }}
\def\sQQ{{\kern 0.27em \vrule height1.45ex width0.03em depth0em
          \kern-0.30em \rm Q}}
\def\qu{{\mathchoice
    {\sQQ}
    {\sQQ}
  {\kern 0.225em \vrule height1.05ex width0.025em depth0em \kern-0.25em \rm Q}
  {\kern 0.180em \vrule height0.78ex width0.020em depth0em \kern-0.20em \rm Q}
        }}
\def\sCC{{\kern 0.27em \vrule height1.45ex width0.03em depth0em
          \kern-0.30em \rm C}}
\def\complex{{\mathchoice
    {\sCC}
    {\sCC}
  {\kern 0.225em \vrule height1.05ex width0.025em depth0em \kern-0.25em \rm C}
  {\kern 0.180em \vrule height0.78ex width0.020em depth0em \kern-0.20em \rm C}
        }}

\newcommand{\ba}{\begin{array}}
\newcommand{\ea}{\end{array}}
\newcommand{\be}{\begin{equation}}
\newcommand{\ee}{\end{equation}}
\newcommand{\bea}{\begin{eqnarray}}
\newcommand{\eea}{\end{eqnarray}}
\newcommand{\beaa}{\begin{eqnarray*}}
\newcommand{\eeaa}{\end{eqnarray*}}

\newcommand{\eps}{\varepsilon}

\def\b{\beta}

\def\z{\zeta}

\font\tenmath=msbm10 \font\sevenmath=msbm7 \font\fivemath=msbm5
\newfam\mathfam \textfont\mathfam=\tenmath
\scriptfont\mathfam=\sevenmath \scriptscriptfont\mathfam=\fivemath

\def \b{\noindent}

\def \={{\buildrel {\rm (law)} \over =}}

\def\qed{ \hfill \vrule width.25cm height.25cm depth0cm\smallskip}

\newcommand{\basa}{\begin{assumption}}
\newcommand{\easa}{\end{assumption}}

\newcommand{\bas}{\begin{assum}}
\newcommand{\eas}{\end{assum}}

\newcommand{\ignore}[1]{}
\begin{document}

\renewcommand{\thefootnote}{\fnsymbol{footnote}}

\renewcommand{\thefootnote}{\fnsymbol{footnote}}

\title{ Absolute continuity of the solution to the stochastic Burgers equation  }
\author{Christian Olivera $^{1,}$\footnote{ C. O. is supported by FAPESP
		by the grant 2020/04426-6	, by  CNPq by the grant
		426747/2018-6 and CAPES by  the grant  
88887.198637/2018-00- MATHAMSUD project SARC (19-MATH-06).} \hskip0.2cm 
Ciprian A. Tudor $^{2}$ \footnote{C. T. acknowledges partial support from the Labex CEMPI (ANR-11-
LABX-0007-01) and MATHAMSUD project SARC (19-MATH-06).}\vspace*{0.1in} \\
$^{1}$ Departamento de Matem\'atica, Universidade Estadual de Campinas,\\
13.081-970-Campinas-SP-Brazil. \\
colivera@ime.unicamp.br \vspace*{0.1in} \\
 $^{2}$  Universit\'e de Lille, CNRS\\
Laboratoire Paul Painlev\'e UMR 8524\\
 F-59655 Villeneuve d'Ascq, France.\\
 \quad tudor@math.univ-lille1.fr\vspace*{0.1in}}
\maketitle

\begin{abstract}
We prove the existence and the Besov regularity of the density of the solution to a general parabolic SPDE which includes the stochastic Burgers equation on an unbounded domain. We use an elementary approach based on the fractional integration by parts. 
\end{abstract}

\medskip

{\bf MSC 2010\/}: Primary 60H15: Secondary 60H05.

 \smallskip

{\bf Key Words and Phrases}:  stochastic Burgersr equation, space-time white noise,    existence  of the density, Besov spaces, fractional integration by parts.

\section{Introduction}

We consider the following parabolic stochastic partial  differential equation (SPDE in the sequel)
\begin{equation}
\label{burgers}
\frac{\partial}{\partial t} u(t,x)= \frac{\partial ^{2}}{\partial x^{2}} u(t, x) -\frac{\partial}{\partial x} g(u(t,x)) + \sigma (t, x, u(t, x)) \frac{\partial ^{2}W}{\partial t \partial x} 
\end{equation}
where the random noise $W$ is a Brownian sheet and the initial condition $u_{0}$ is deterministic. Our purpose is to cover both the situations  when  the space variable belongs  to the whole real line and to  the unit interval $[0,1]$.  That is, will assume in (\ref{burgers}) that 
$$(t,x) \in [0, T] \times I \mbox{ with } I= \mathbb{R} \mbox{ or } I=[0,1].$$
In both situations, we will assume 
$$u_{0}(t,x)= u_{0}(x) \mbox{ for every } t\in [0, T], x\in I $$
while for $I=[0,1]$  we also impose Dirichlet boundary conditions 
$$u_{0}(t, 0)= u_{0}(t,1) =0 \mbox{ for every } t\in [0, T]. $$  

In the particular case $g(x)=\frac{1}{2}x^{2}$, the equation (\ref{burgers}) is called {\it the  stochastic  Burgers equation } and it has been studied by many authors. It is one of the most known singular stochastic partial differential equations.  For  a complete description of the   deterministic Burgers equation (i.e. $\sigma\equiv0$ in (\ref{burgers})) we refer to the monograph \cite{Bur} and the references therein.  The Burgers equation with a random forcing term has    been also studied by several authors, see e.g. \cite{AS}, \cite{BCJ}, \cite{BGN}, \cite{CO},    \cite{dPDT},   \cite{Dermoune},
\cite{Gy}, \cite{Gu},  \cite{HLOUZ1}, \cite{Kim},   \cite{LNP}, \cite{LeNu}, \cite{Mo}. The interest to study such stochastic equations  comes from  the fact during the past few decades, the stochastic Burgers equation has found applications in various  fields ranging from statistical physics, cosmology, and fluid dynamics to engineering. In particular, the problem of  Burgers turbulence, that is,  the study of the solutions to the Burgers equation with random  initial conditions or with a  random noise  is a central issue in the study of nonlinear systems out of equilibrium (see  e.g. the surveys \cite{Beck} and \cite{Wie}).

Let us recall some known facts concerning the existence, the uniqueness and the properties of the solution to (\ref{burgers}). For instance, in the case of the unbounded spatial domain, the existence and the pathwise properties of the solution have been studied  in \cite{BCJ}  (for the additive noise, i.e. $\sigma \equiv 1$),  \cite{GyNu}, \cite{LNP}, while for the case of a bounded domain we refer to \cite{LNP} (for the particular case of the stochastic Burgers equation  $g(x)= \frac{1}{2}x^{2}$) or \cite{Gy}, \cite{Mo} for a more general situation. 

Our purpose is to discuss the absolute continuity of the law of the solution to (\ref{burgers}) with a particular focus on the case of the Burgers equation. The existence and some regularity of the density of the  solution are  already known in the case of a bounded domain. Indeed, such results have been obtained in \cite{LNP} or \cite{Mo} by using the techniques of the Malliavin calculus, which requires rather stronger assumptions on the coefficients and technical proofs. As far as we know, there are no results concerning the absolute continuity of the law of the solution  when the space variable belongs to the whole real line. 

In order to prove the existence of a density for the solution to the parabolic SPDE (\ref{burgers}), we propose here a different approach,  based on the fractional integration by parts technique developed in \cite{Ro} or \cite{DeRo}. This method is rather simple, it fits well with some classes of stochastic (partial) differential equations and it allows to obtain some  properties of the density that cannot be obtained via Malliavin calculus or other techniques. Examples of applications of the integration by parts technique to the study of densities of solutions to various types of stochastic equations can be found in \cite{Deb}, \cite{DeRo}, \cite{Fou}, \cite{Ro}, \cite{Sanz1}, \cite{OT1}. We obtain via a simple application of this method the existences and the Besov regularity of the density of the solution to (\ref{burgers}), including the case of the Burgers equation $g(x)= \frac{1}{2}x^{2}$. As we pointed out before, such a result seems to be completely new  for the unbounded domain while it provides a simpler alternative to  the more technical proofs  in  \cite{LNP} or \cite{Mo} in the case $I=[0, 1]$.  

We organized our work as follows. Section 2 includes some general facts on the parabolic SPDEs of type (\ref{burgers}), on its associated fundamental solution (Green kernel) and the integration by parts (or "smoothing"
 lemma). In Section 3 we apply the fractional integration by parts technique to get the  existence and the regularity in Besov dense for the density of the solution to (\ref{burgers}) while in Section 4 we discuss several examples, including the case of the Burgers equation on the real line.

By $c, C$ we denote generic strictly positive constants that are allowed to change from one line to another (and even on the same line). 

\section{Preliminaries}
In this preliminary part, we define the noise and the solution to the parabolic SPDE (\ref{burgers}). We also recall the fractionla integration by parts lemma which is the key element to obtain the existence and the Besov regularity of the density.

\subsection{The equation}

Let $I=[0,1]$ or $I=\mathbb{R}$ and let  $ \left(W(t,x), t\in [0,T], x\in I\right)$ be a centered Gaussian field defined on a complete probability space $\left(\Omega, \mathcal{F}, P\right)$ with covariance, for every $s,t\in [0, T]$ and $x,y\in I$
\begin{equation*}
\mathbf{E} W(s,x) W(t, y)= (s\wedge t) (\vert x\vert \wedge \vert y\vert ) 1_{ (0, \infty)}(xy)  \mbox{ if } I=\mathbb{R} 
\end{equation*}
and 
\begin{equation*}
\mathbf{E} W(s,x) W(t, y)= (s\wedge t) (x\wedge y)  \mbox{ if } I=[0,1].
\end{equation*}
We will denote by $\mathcal{F}_{t}$ the sigma-algebra generated by the random variables $\left( W(s, x), 0\leq s\leq t, x\in I\right) $. We assume that $\mathcal{F}= \mathcal{F}_{T}$ and we will say that a process $ \left( Y(t,x) , t\in [0,T], x\in I\right) $ is adapted with respect to the filtration $\left( \mathcal{F}_{t}\right) _{t\in [0, T]}$ if the random variable $Y(t,x) $ is $\mathcal{F}_{t}$-measurable for every $t\in [0, T], x\in I$.

The solution to (\ref{burgers}) is understood in the mild sense, that is $u$ is a $\mathcal{F}_{t}$-adapted stochastic  field that  satisfies the integral equation 

\begin{eqnarray}
u(t,x)&=& \int_{I} G_{t}(x-y) u_{0} (y)dy + \int_{0} ^{t} \int_{I} \frac{\partial }{\partial y} G_{t-s}(x-y) g(u(s,y)) dy ds \nonumber\\
&&+ \int_{0} ^{t} \int_{I} G_{t-s}(x-y) \sigma ( u(y, u(s, y))) W (ds, dy), \hskip0.4cm (t,x)\in [0, T]\times I.   \label{mild}
\end{eqnarray} 
where the stochastic integral $W(ds, dy)$ is the Walsh integral with respect to Brownian sheet $W$ (see \cite{Walsh} for its definition and properties) and the Green kernel $G$ is the fundamental solution associated to the heat  equation. It has the following expression (upon the situations $I=\mathbb{R}$ or $I=[0,1]$): for $t>0$,
\begin{equation}
\label{G1}
G_{t}(x):=G_{1, t} (x)=\frac{1}{\sqrt{4\pi t}} e ^{-\frac{x^{2}}{4t}} \mbox{ if } I=\mathbb{R}
\end{equation}
and
\begin{equation}
\label{G2}
G_{t}(x):=G_{2, t}(x)=\frac{1}{\sqrt{4\pi t}} \sum_{n=-\infty} ^{\infty} \left[ e ^{-\frac{(y-x-2n)^{2}}{4t}}- e ^{-\frac{(y+x-2n)^{2}}{4t}} \right] \mbox{ if } I=[0, 1].
\end{equation}

The following properties of the Green kernel will be needed later.

\begin{lemma}\label{ll1}
Let $G_{1}, G_{2}$ be given by (\ref{G1}), (\ref{G2}) respectively. Then we have
\begin{enumerate}
\item  For  every $0<\eps <t$, we have
\begin{equation*}
 \int_{t-\eps}  ^{t}  \int_{\mathbb{R}} G_{1, t-s} (x-y) dyds \geq C_{2} \sqrt{\eps} \mbox{ and }
\int_{t-\eps}  ^{t} \int_{0}^{1} G_{2, t-s} (x-y) dyds \geq C \sqrt{\eps}.
\end{equation*}

\item For every $x\in I$ and $t\in (0, T]$, 
\begin{equation*}
\int_{\mathbb{R}} G^ {2}_{1, t} (x-y) dy \leq C t ^{-\frac{1}{2}} \mbox{ and } \int_{\mathbb{R}} G^ {2}_{2, t} (x-y) dy \leq C t ^{-\frac{1}{2}}.
\end{equation*}

\item 
For every $x\in I$ and $t\in (0, T]$, 
\begin{equation*}
  \int_{\mathbb{R}} \left| \frac{ \partial }{\partial y} G_{1, t}(x-y) \right| dy  \leq C t ^{-\frac{1}{2}} \mbox{ and  } 
 \int_{0}^{1} \left| \frac{ \partial }{\partial y} G_{2, t}(x-y) \right| dy  \leq C t^{-\frac{1}{2}}.
\end{equation*}
\end{enumerate}
\end{lemma}
{\bf Proof: }The  proof of the lower bound of  point 1.   is immediate for $I=\mathbb{R}$, since by (\ref{G1})

\begin{eqnarray}
&&\int_{t-\eps} ^{t}\int_{I} G^{2}_{1, t-s}(x-y) dyds= C \int_{t-\eps}^{t} (t-s) ^{-1} \int_{\mathbb{R}} e ^{-\frac{(x-y) ^{2} }{2(t-s)} }dyds\nonumber\\
&=&C \int_{t-\eps}^{t} (t-s)^{-\frac{1}{2} }ds= C \eps^{\frac{1}{2}}\label{10f-5}
\end{eqnarray}
while for $I=[0,1] $ and $G=G_{2}$, the inequality in the statement  has been proved in  Lemma A.2 in \cite{Mo}. Concerning point 2.,  by relation (2.4) in \cite{LNP}, with some $a, b>0$,
$$\left| G_{2, t}(x, y) \right|\leq at ^{-\frac{1}{2}} e ^{-\frac{b\vert x-y\vert ^{2}}{t}}$$
and  then 
$$ \int_{0}^ {1} G_{2, t}^{2} (x-y) dy \leq a^ {2} t^{-1}\int_{\mathbb{R}}e ^{-\frac{2b\vert x-y\vert ^{2}}{t}}dy\leq Ct ^{-\frac{1}{2}}.$$
The same calculations are obvious for $G_{1}$, since
$$ \int_{\mathbb{R}} G_{1, t}^{2} (x-y) dy = C t ^ {-1} \int_{\mathbb{R}}e ^{-\frac{\vert x-y\vert ^{2}}{2t}}dy \leq Ct ^{-\frac{1}{2}}.$$

For point 3., when $I=[0,1] $ we have by relation (2.5) in \cite{LNP}, for some $a,b>0$, abnd for every $x,y \in I$, $ t\in (0, T]$, 

$$ \left| \frac{\partial}{\partial y} G_{2, t-s}(x-y)\right|\leq at ^{-1} e ^{-\frac{b\vert x-y\vert ^{2}}{t}}$$
and thus
$$\int_{0}^{1} \left| \frac{ \partial }{\partial y} G_{2, t-s}(x-y) \right| dy \leq C t ^{-1}\int_{\mathbb{R}}e ^{-\frac{b\vert x-y\vert ^{2}}{t}}dy\leq Ct ^{-\frac{1}{2}}.$$

\b  When $I=\mathbb{R}$, the above two relations also also true for the kernel $G_{1}$.

\qed

The existence and the uniqueness of the solution has been obtained in e.g. \cite{BCJ}, \cite{Gy}, \cite{GyNu}.  The conditions for the existence and uniqueness are slightly different in the cases of the bounded or unbounded domain. In order to kep an unitary approach, we will first show  the absolute continuity of the law of (\ref{mild}) by assuming the existence and the pathwise regularity  of the solution. In the last part of the paper which contains the examples, we will describe for each case the conditions  that ensure existence, uniqueness and other properties of the solution.

\subsection{Besov spaces and the fractional integration by parts}

The monograph  \cite{Triebel}  offers  a complete exposition on Besov spaces. In this work we only need to consider the particular  Besov space $\mathcal{B} _{1, \infty}^{s} $ with $s>0$, which is defined below.

Consider a function $f:\mathbb{R} \to \mathbb{R}$ and for every $x,h\in \mathbb{R}$, put
$ (\Delta _{h} ^{1} f) (x)= f(x+h) -f(x)$
 and for $n\geq 1$ integer, define the $n$th increment of the function $f$ at lag $h$ by 

\[
(\Delta_{h}^{n}f )(x)= \Delta_{h}^{1}\left(\Delta_{h}^{n-1}f\right))(x)= \sum_{j=0}^{n} (-1)^{n-j} \binom{n}{j} f(x+jh). 
\]

\b  For instance, if $m=2$, $(\Delta_{h}^{n}f )(x)=f(x+2h)-2f(x+h)+ f(x)$.  For $0<s<n$ we define the norm
\begin{equation}\label{besovn}
\Vert f\Vert _{\mathcal{B} _{1,\infty}^{s}}=\Vert f\Vert _{L ^{1}(\mathbb{R} )}+ \sup_{\vert h\vert \leq 1} \vert h\vert ^{-s} \Vert \Delta _{h} ^{n} f\Vert _{ L ^{1}(\mathbb{R} ^{d})}.
\end{equation}
It can be shown that for any $n,m>s$, the norms obtained in (\ref{besovn}) using different  $n$ and $m$ are equivalent. Therefore, one can define the Besov space $\mathcal{B} _{1, \infty}^{s} $ as the set of functions $f\in L ^{1} (\mathbb{R} )$ such that 
$\Vert f\Vert _{\mathcal{B} _{1,\infty}^{s}}<\infty.$

For $\alpha \in (0, 1)$, let $C_{b} ^{\alpha} (I)$ denote the set of bounded  measurable functions $h: I \to \mathbb{R}$ such that 
\begin{equation*}
\Vert h\Vert _{C_{b} ^{\alpha} (I)} := \Vert h\Vert _{\infty}+ \sup_{x, y\in I, x\not=y, \vert x-y\vert \leq 1}\frac{ \vert h(x)- h(y)\vert}{\vert x-y\vert ^{\alpha}} <\infty.
\end{equation*}

Our main tool to get the existence and the regularity of the density of the solution to (\ref{burgers}) is the following smoothing lemma  from \cite{Ro} (also called {\it the fractional integration  by parts lemma}).

\begin{lemma}\label{Ro}
 Let $X$ be a  real valued random variable. If there exist
an integer $m\geq 1$, two  real numbers $s>0, \gamma>0$, with 
$\gamma <s <m$, and a constant $K>0$ such that for every $\phi\in C_{b}^{\gamma}(\mathbb{R})$ and 
$h\in \mathbb{R}$, with $|h|\leq 1$,

\[
\mathbf{E}\left[\Delta_{h}^{m}\phi(X)\right]\leq K |h|^{s} \|\phi\|_{C_{b}^{\gamma}(\mathbb{R})}, 
\]
then $X$ admits a density $f_{X}$ with respect to Lebesgue measure on $\mathbb{R}$. Moreover, $f_{X}\in B_{1,\infty}^{s-\alpha}$
and $
\| f\|_{B_{1,\infty}^{s-\alpha}}\leq C (1+ K). $

\end{lemma}

\section{The existence  and Besov regularity of the density}

We will use the fractional integration by parts (Lemma \ref{Ro}) in order to obtain the absolute continuity of the law of the random variable $u(t,x)$ given by (\ref{mild}), for fixed $t\in (0, T]$ and $x\in I$. We will state and prove a general result, that includes both the cases of bounded and unbounded domains, by assuming the existence, the H\"older regularity in time  and some moment estimates for the solution to (\ref{burgers}). In the next section, we will treat separately several particular cases and for each of them we will show that these properties are satisfied.

Let us consider the following assumptions: there exist $L_{1}, L_{2}, K>0$ such that 
\begin{itemize}
\item 
The diffusion coefficient $\sigma:  I \times \mathbb{R} \to\mathbb{R}$  does not depend on time and it  is globally Lipschitz i.e.
\begin{equation}\label{lip}
\left| \sigma (x, r)- \sigma (x, v) \right| \leq L_{1} \vert r-v\vert \mbox{ for every } r, v\in \mathbb{R}, x\in I.
\end{equation}

and 
\begin{equation}
\label{cc1}\vert  \sigma (x, r)\vert ^{2} \geq K \mbox{ for every } x\in I, r\in\mathbb{R}.
\end{equation}

\item The function $g:\mathbb{R}\to \mathbb{R}$ satisfies 

\begin{equation}\label{cc2}
\vert g(r)-g(v)\vert\leq L_{2} (1 +\vert r\vert +\vert v\vert ) \vert r-v\vert \mbox{ for every } r, v\in \mathbb{R}. 
\end{equation}

\end{itemize}

We assume for simplicity that $\sigma$ does not depend on time. The case of a time-depedent diffusion coefficient can be treated in a similar way, by assuming  the Lipschitz condition in time for it. Notice that the above conditions (\ref{lip}), (\ref{cc1}) and (\ref{cc2}) does not guarantee the existence of the solution (\ref{mild}), see Section \ref{sec4}.

Let us state our general result concerning the absolute continuity of the law of the solution to (\ref{burgers}).
\begin{theorem}
\label{tt1} Assume (\ref{lip}), (\ref{cc1}) and (\ref{cc2}). Also assume that the equation (\ref{burgers}) admits a unique solution $\left( u(t,x), t\in [0, T], x\in I\right)$  and there exists $\beta >0$ such that 
for every $p\geq 2$, 
\begin{equation}
\label{holder}\mathbf{E} \left| u(t,x)- u(s, x)\right| ^{p} \leq c \vert t-s\vert ^{\beta p }
\end{equation}
for every $x \in I$, $t,s\in [0, T]$ and 
\begin{equation}
\label{momp}\sup_{t\in [0, T], x\in I} \mathbf{E}\vert u(t, x) \vert ^{p}\leq C_{T}. 
\end{equation}
 Then  for every $t\in (0, T]$ and $x\in I$, the random variable $u(t,x)$ admits a density which belongs to the Besov space $\mathcal{B}_{1, \infty}^{p}$ with $p<2\beta$.

\end{theorem}
{\bf Proof: } Fix $t\in (0, T]$ and $0<\eps<\frac{t}{2}$. We introduce the auxiliary processes
\begin{eqnarray}
u_{0, \eps}(t,x)&=& \int_{I} G_{t}(x-y) u_{0} (y)dy + \int_{0} ^{t-\eps} \int_{\mathbb{I}} \frac{\partial }{\partial y} G_{t-s}(x-y) g(u(s,y)) dy ds \nonumber  \\
&&+ \int_{0} ^{t-\eps} \int_{I} G_{t-s}(x-y) \sigma ( u(y, u(s, y))) W (ds, dy)   \label{aux1}
\end{eqnarray} 
and 

\begin{eqnarray}
u_{\eps} (t,x) &=& u_{0, \eps } (t,x)+ \int_{t-\eps} ^{t} \int_{I} \frac{\partial }{\partial y} G_{t-s}(x-y) g(u(t-\eps, y))dyds \nonumber  \\
&&+  \int_{t-\eps} ^{t} \int_{I} G_{t-s}(x-y) \sigma(u(y, u(t-\eps), y))W (ds, dy).\label{aux2}
\end{eqnarray}

In order to apply Lemma \ref{Ro},  we need to calculate 
$$ \mathbf{E} (\Delta _{h}^{m} \varphi ) (u(t,x))$$
for every $h>0$, $m\geq 1$ integer and $\varphi: \mathbb{R}\to \mathbb{R} $ a function  in $C_{b} ^{\gamma }(\mathbb{R})$   for $\gamma \in (0, 1)$. We write 

\begin{eqnarray}
\mathbf{E} (\Delta _{h}^{m} \varphi ) (u(t,x))&=& \mathbf{E} (\Delta _{h}^{m} \varphi ) (u_{\eps} (t,x)) + \mathbf{E} \left[  (\Delta _{h}^{m} \varphi ) (u(t,x))- (\Delta _{h}^{m} \varphi ) (u_{\eps}(t,x))\right]\nonumber \\
&=:&(PE) _{\eps}(m,h,t,x)+  (AE) _{\eps}(m,h,t,x)\label{8f-1}
\end{eqnarray}
with $u_{\eps}$ given by (\ref{aux2}). The first summand above is usually called {\it the probabilistic estimate } while the second term is the {\it approximation error. } Let us estimate separately these two quantities.

\vskip0.3cm

\b {\bf  Calculation of the probabilistic estimate. } Only condition (\ref{cc1}) is needed for this estimate.  We write, for every $t\in (0,T], x\in I$ and $0<\eps <\frac{t}{2}$, 
$$u_{\eps} (t,x) = Z_{t, x, \eps}+ I_{t,x, \eps}$$
with
\begin{equation}
\label{ze}
Z_{t, x,\eps} = u_{0, \eps } (t,x)+ \int_{t-\eps} ^{t} \int_{I} \frac{\partial }{\partial y} G_{t-s}(x-y) g(u(t-\eps, y))dyds 
\end{equation}
and
\begin{equation}
\label{ie}
I_{t,x, \eps}= \int_{t-\eps} ^{t} \int_{I} G_{t-s}(x-y) \sigma(u(y, u(t-\eps), y))W (ds, dy).
\end{equation}
From (\ref{8f-1})
\begin{eqnarray*}
(PE) _{\eps}(m,h,t,x)&=& \mathbf{E} (\Delta _{h}^{m} \varphi ) (u_{\eps} (t,x)) \\
&=& \mathbf{E}\left[  \mathbf{E}\left( (\Delta _{h}^{m} \varphi )(u_{\eps} (t,x)) /\mathcal{F}_{t-\eps} \right)\right] = \mathbf{E}\left[  \mathbf{E}\left( (\Delta _{h}^{m} \varphi )(Z_{t, x, \eps}+ I_{t,x, \eps}) /\mathcal{F}_{t-\eps} \right)\right].
\end{eqnarray*}

Notice that the random variable $Z_{t,x, \eps}$ given by (\ref{ze}) is measurable with respect to $\mathcal{F}_{t-\eps} $ while $I_{t,x, \eps}$ given by (\ref{ie}) is,  conditionally  on $\mathcal{F}_{t-\eps} $, independent of $Z_{t,x, \eps}$. Thus

\begin{equation}
\label{10f-1}(PE) _{\eps}(m,h,t,x)= \mathbf{E} f(Z_{t,x, \eps})
\end{equation}
with
$$f(y)= \mathbf{E}\left[ (\Delta _{h}^{m} \varphi)(y+ I_{t, x, \eps})\right] \mbox{ for every } y\in \mathbb{R}. $$

Also, we can se that the random variable  $I_{t,x, \eps}$ is, conditionally on $\mathcal{F}_{t-\eps} $, a Gaussian random variable with zero expectation and variance 
\begin{equation*}
V_{t,x, \eps}= \int_{t-\eps}^{t} \int_{I} G^{2}_{t-s}(x-y)\sigma ^{2} (y, u(t-\eps, y))dyds
\end{equation*} So the conditional density of $I_{t,x, \eps}$ is
$$g_{t,x ,\eps} (z)= \frac{1}{\sqrt{2\pi V_{t,x, \eps}}}e ^{-\frac{z^{2}}{{2V_{t,x, \eps}}}}, \hskip0.5cm z\in \mathbb{R}.$$

\b Therefore
$$f(y)=\int_{\mathbb{R}} (\Delta _{h}^ {m} \varphi)(y+z) g_{t,x, \eps}(z)dz= \int_{\mathbb{R}} \varphi (y+z) \left( \Delta _{-h} ^{m} g_{t, \eps }(z) \right) dz 
\leq \Vert \varphi\Vert _{\infty} \Vert \Delta _{-h} ^{m} g_{t, \eps }(z) \Vert _{L ^{1} (\mathbb{R} ) }.
$$

We know  from  \cite{Ro} (see page 5, two lines after (2.7)) that 
\begin{equation}\label{18m-2}
 \Vert \Delta _{-h} ^{m} g_{t, \eps (x) }\Vert _{L ^{1} (\mathbb{R} ) }\leq C\left(  \frac{\vert h\vert }{\sqrt{V_{t,x, \eps}} }\right) ^{m}  \leq   C \left(  \frac{\vert h\vert }{\eps ^{\frac{1}{2}}}\right) ^{m} 
\end{equation}
for every $h>0$ and for any integer $m\geq 1$. For the last inequality in (\ref{18m-2}) we used the fact that by  assumption (\ref{cc1}) and  Lemma \ref{ll1} point 1., for every $t,x, \eps$,
$$V_{t,x,\eps}\geq K \int_{t-\eps} ^{t}\int_{I} G^{2}_{t-s}(x-y) dyds \geq C \sqrt{\eps}. $$
 Thus 
\begin{equation}
\label{10f-2}
\sup_{y\in \mathbb{R}} \vert f(y) \vert \leq C \left(  \frac{\vert h\vert }{\eps ^{\frac{1}{2}}}\right) ^{m}.
\end{equation}
and (\ref{10f-1}) and (\ref{10f-2}) imply for every $h>0$,  $m\geq 1$ integer  and $(t,x)\in (0, T]\times I$,
\begin{equation}
\label{10f-3}
\left| (PE) _{\eps}(m,h,t,x)\right| \leq C \left(  \frac{\vert h\vert }{\eps ^{\frac{1}{2}}}\right) ^{m}.
\end{equation}

\qed

\vskip0.3cm

\noindent {\bf Calculation of the approximation error. } To handle the term $(AE) _{\eps}(m,h,t,x)$ in (\ref{8f-1}), we will  use the assumptions  (\ref{lip}), (\ref{cc2}), (\ref{holder}) and (\ref{momp}). Recall that by (\ref{8f-1})

$$(AE) _{\eps}(m,h,t,x)= \mathbf{E} \left[  (\Delta _{h}^{m} \varphi ) (u(t,x))- (\Delta _{h}^{m} \varphi ) (u_{\eps}(t,x))\right].$$

\b From (\ref{mild}) and (\ref{aux2}),  for every $0<\eps <\frac{t}{2}$ and for every $x\in I$, 

\begin{eqnarray*}
u(t,x)- u_{\eps}(t,x)&=& \int_{t-\eps}^{t} \int_{I} \frac{\partial}{\partial y} G_{t-s}(x-y) \left( g(u(s,y))- g(u(t-\eps, y)) \right) dyds \\
&&+ \int_{t-\eps}^{t} \int_{I} G_{t-s}(x-y) \left( \sigma (u(y, u(s,y))- \sigma (y, u(t-\eps, y))\right) W(ds, dy).
\end{eqnarray*}
This implies, since $\gamma \in (0,1)$, 
\begin{equation}
\label{10f-4} \left| (AE) _{\eps}(m,h,t,x)\right| \leq C\Vert \varphi \Vert _{ C_{b} ^{\gamma}(I) }   \mathbf{E} \vert u(t,x)- u_{\eps}(t,x)\vert ^{\gamma}\leq C  \Vert \varphi \Vert _{ C_{b} ^{\gamma} (I)}  \left( \mathbf{E}\vert u(t,x)- u_{\eps}(t,x)\vert \right) ^{\gamma}.
\end{equation}

We then need to estimate $\mathbf{E}\vert u(t,x)- u_{\eps}(t,x)\vert $.  We can write,

\begin{eqnarray*}
\mathbf{E}\vert u(t,x)- u_{\eps}(t,x)\vert &\leq&  \int_{t-\eps}^{t} \int_{I} \left| \frac{\partial}{\partial y} G_{t-s}(x-y)\right| \mathbf{E}\left| g(u(s,y))- g(u(t-\eps, y)) \right| dyds\\
&&+\mathbf{E}\left|  \int_{t-\eps}^{t} \int_{I} G_{t-s}(x-y) \left( \sigma (u(y, u(s,y))- \sigma (y, u(t-\eps, y))\right) W(ds, dy)\right| \\
&:=& A+B.
\end{eqnarray*} 

Now, by (\ref{cc2}), (\ref{momp})  and H\"older's inequality

\begin{eqnarray*}
A &\leq& C \int_{t-\eps}^{t} \int_{I} \left| \frac{\partial}{\partial y} G_{t-s}(x-y)\right| \mathbf{E} (1+ \vert u(s,y)\vert + \vert u(t-\eps, y)\vert ) \vert u(s, y) - u(t-\eps, y) \vert   dyds\\
&\leq& C \int_{t-\eps}^{t} \int_{I}  \left| \frac{\partial}{\partial y} G_{t-s}(x-y)\right|\left(  \mathbf{E} \vert u(s, y) - u(t-\eps, y) \vert ^{2} \right) ^{\frac{1}{2}}
\end{eqnarray*}
and by using the H\"older assumption (\ref{holder}) with $p=2$, we obtain 

\begin{eqnarray*}
A &\leq & C \int_{t-\eps}^{t} \int_{I}  \left| \frac{\partial}{\partial y} G_{t-s}(x-y)\right|(s-t+\eps) ^{\beta}dyds \leq  C \eps ^{\beta }  \int_{t-\eps}^{t} \int_{I}  \left| \frac{\partial}{\partial y} G_{t-s}(x-y)\right|dyds.
\end{eqnarray*}
By Lemma \ref{ll1} point3. and using $0<\eps <\frac{t}{2}$ with $t\in (0, T]$, 

\begin{equation}
\label{10f-8}A \leq C\eps ^{\beta } \int_{t-\eps}^{t} (t-s) ^{-\frac{1}{2}} ds \leq C\eps ^{\beta } \int_{t-\eps}^{t} s ^{-\frac{1}{2}} ds= C \eps ^{\beta} \left( t ^{\frac{1}{2}} - (t-\eps) ^{\frac{1}{2}}\right)\leq C t^{-\frac{1}{2}} \eps ^{\beta +1}.
\end{equation}
 
\b For the term $B$, we have by (\ref{lip}) and (\ref{holder})
\begin{eqnarray*}
B ^{2} &\leq& 
\int_{t-\eps}^{t} \int_{I} G^{2}_{t-s}(x-y) \mathbf{E}\left( \sigma (u(y, u(s,y))- \sigma (y, u(t-\eps, y))\right)^{2} dyds\\
&\leq & C\int_{t-\eps}^{t} \int_{I} G^{2}_{t-s}(x-y)\mathbf{E} \left| u(s,y)-u(t-\eps, y)\right| ^{2} dyds \\
&\leq & C \int_{t-\eps}^{t} \int_{I} G^{2}_{t-s}(x-y) (s-t+\eps) ^{2\beta } dyds \\
&\leq & C \eps ^{2\beta } \int_{t-\eps}^{t}  (t-s) ^{-\frac{1}{2}}dyds \leq Ct ^{-\frac{1}{2}}\eps ^{2\beta  +1} 
\end{eqnarray*}
where we used Lemma \ref{ll1} point 2. and we  proceeded as for the bound (\ref{10f-8}).  Therefore, by (\ref{10f-4})
\begin{equation}
\label{10f-6}\left| (AE) _{\eps}(m,h,t,x)\right| \leq C\Vert \varphi \Vert _{ C_{b} ^{\gamma }(I) }   \eps ^{\gamma \left( \beta +\frac{1}{2}\right)}.
\end{equation}
We obtained by (\ref{10f-3}) and (\ref{10f-6}), for every $h>0$, $m\geq 1$ integer and $(t,x)\in [0, T] \times I$, 

\begin{equation}
\label{14f-1}
\vert \mathbf{E} (\Delta _{h}^{m} \varphi ) (u(t,x))\vert \leq C\Vert \varphi \Vert _{ C_{b} ^{\gamma } (I)}   \left[ \eps ^{ \left( \beta +\frac{1}{2}\right)}+\left(  \frac{\vert h\vert }{\eps ^{\frac{1}{2}}}\right) ^{m} \right]^ {\gamma} .
\end{equation}

We choose $\eps = h ^{A} $ with suitable $A$ which optimizes  the right-hand side  of (\ref{14f-1}). This gives
$$A= \frac{2m}{m+2 (\beta + \frac{1}{2})}.$$
and so  (\ref{14f-1}) becomes
$$\vert \mathbf{E} (\Delta _{h}^{m} \varphi ) (u(t,x))\vert \leq C\Vert \varphi \Vert _{ C_{b} ^{\alpha} (I)} h ^{s}  \mbox{ with }s= \frac{2m\gamma (\beta +\frac{1}{2})}{m+2 (\beta +\frac{1}{2})}.$$

When $m$ is large enough,  we see that $s$ becomes arbirarly close to $2\gamma\left( \beta +\frac{1}{2}\right)$.  Then, by Lemma \ref{Ro}, the solution to (\ref{burgers}) admits a density and  this  density belongs to the Besov space $\mathcal{B} _{1, \infty} ^{p}$ with 
$$p< s-\gamma =2\gamma \left( \beta +\frac{1}{2}\right)- \gamma= 2\gamma \beta.$$
By choosing $\gamma $ arbitrarly close to $1$, we get that the density belongs to the Besov space $\mathcal{B} _{1, \infty} ^{p}$ with $0<p<2\beta  $. \qed

Notice that the Besov  regularity of the density is related to the H\"older regularity in time of the solution. The  spatial regularity (which is also known from \cite{Mo} or \cite{LeNu}) of the solution does not affect the regularity of the density.

\section{Applications}\label{sec4}
Let us now discuss the stochastic equation (\ref{burgers}) in some particular cases. 

\subsection{The stochastic Burgers equation on $\mathbb{R}$}

We will show that our result stated in Theorem \ref{tt1} can be applied to the stochastic Burgers equation on the whole real line. Recall  that no results are known in this case on the existence of the density of the solution. 

We will assume throughout this section that  
\begin{equation}
\label{cc4}
I=\mathbb{R} \mbox{ and } g(x)= \frac{1}{2}x ^ {2}  \mbox{ for every } x \in \mathbb{R}
\end{equation}
which corresponds to the case of the stochastic Burgers equation. 

In order to ensure the existence of the solution (\ref{mild})  in the situation (\ref{cc4}), let us consider the following additional conditions: 
\begin{itemize}
\item For some $\alpha \in (0, 1)$, the initial condition $u_{0} $ satisfies

\begin{equation}\label{u0}
u_{0} \in L^{1} (\mathbb{R})\cap L^{2} (\mathbb{R}) \cap C_{b} ^{\alpha} (\mathbb{R}).
\end{equation}

\item  The diffusion coefficient $\sigma$ is such that 
\begin{equation}\label{cc3}
\vert \sigma (x, r) \vert \leq f(x) \mbox{ for every } x,r \in \mathbb{R}
\end{equation}
where $f$ is a non-negative function such that $f\in L^ {2} (\mathbb{R})\cap L^ {q}(\mathbb{R}) $ with $q>2$.

\end{itemize}

The following result has been proved in  \cite{GyNu} (the existence of the solution) and in \cite{LNP} (the H\"older regularity and the moment estimates of the solution). 

\begin{theorem}Assume (\ref{lip}),  (\ref{cc4}),  (\ref{u0}) and  (\ref{cc3}). Then (\ref{burgers}) admits an unique solution which satisfies, for all $s, t \in [0, T]$,  $x, y\in \mathbb{R}$  and $p\geq 2$

\begin{equation*}
\mathbf{E} \left| u(t,x)- u(s, y) \right|  ^ {p} \leq C \left( \vert t-s\vert ^ {\frac{\alpha}{2} \wedge (\frac{1}{4}-\frac{1}{2q})}+ \vert x-y\vert ^ {\alpha \wedge (\frac{1}{2} -\frac{1}{q})}\right)^ {p}
\end{equation*}
and
\begin{equation*}
\sup_{t\in [0, T], x\in \mathbb{R}} \mathbf{E} \vert u(t, x) \vert ^{p}\leq C_{T}.
\end{equation*}

\end{theorem}

We obtain, from the above result and Theorem \ref{tt1} the absolute continuity of the the law of the stochastic Burgers equation on an unbounded domain.

\begin{corollary}
Assume (\ref{cc1}), (\ref{lip}),  (\ref{u0}), (\ref{cc3}) and (\ref{cc4}),  and let $(u(t,x), t\in [0, T], x\in \mathbb{R})$ be the unique solution to (\ref{burgers}). Then for every $t\in (0, T], x\in \mathbb{R}$, the random variable $u(t,x)$ admits a density  (denoted $f_{t, x}$) which satisfies
$$f_{t, x} \in \mathcal{B} _{1, \infty}^ {p} \mbox{ with } p< \alpha \wedge \left(\frac{1}{2}-\frac{1}{q}\right).$$ 
\end{corollary}

\begin{remark}
By choosing a regular enough initial condition (for example, $u_{0}$ is H\"older continuous or order $\alpha > \frac{1}{2}$), we obtain that the density of the solution to the Burgers equation  belongs to the Besov space  $\mathcal{B} _{1, \infty}^ {p} \mbox{ with } p<  \frac{1}{2}-\frac{1}{q}$ with $q$ from   (\ref{cc3}). The H\"older regularity in space of the solution does not affect the regularity of the density. 
\end{remark}

\subsection{The Burgers type equation on a bounded domain}

We will regard the  density of the solution to (\ref{burgers}) when the space variable belongs to the unit interval $[0, 1]$.  We will restrict to the case when the function $g$ in (\ref{burgers}) is Lipschiz continuous. This is because, in order to apply Theorem \ref{tt1} we need the H\"older continuity  of the solution which has been proved, as far as we  know, only for $g$ Lipschitz  (see \cite{Mo}).

Let us consider the following assumptions: 
\begin{itemize} 
 \item For some $\alpha \in (0, 1)$, the initial condition $u_{0} $ satisfies

\begin{equation}\label{u00}
I=[0,1] \mbox{ and } u_{0} \in C_{b} ^{\alpha} (I).
\end{equation}

\item The function $g$ satisfies
\begin{equation}
\label{cc5}
 \vert g(x)-g(y)\vert \leq C \vert x-y\vert \mbox{ for every } x, y\in [0, 1].
\end{equation}

\end{itemize} 

The next result is Theorem 2.1 in \cite{Mo}.

\begin{theorem}
Assume (\ref{lip}), (\ref{cc1}),  (\ref{u00}) and (\ref{cc5}). Then the stochastic equation (\ref{burgers})  admits a unique solution which satisfies, for every $s, t\in [0, T], x,y\in [0,1] $ and $p\geq 2$ 
\begin{equation}\label{10f-9}
\mathbf{E} \left| u(t,x)- u(s, y) \right|  ^ {p} \leq C \left( \vert t-s\vert ^ {\frac{\alpha}{2} \wedge \frac{1}{4}}+ \vert x-y\vert ^ {\alpha \wedge \frac{1}{2} }\right)^ {p}
\end{equation}
and
\begin{equation}\label{10f-10}
\sup_{t\in [0, T], x\in [0,1} \mathbf{E} \vert u(t, x) \vert ^{p}\leq C_{T}.
\end{equation}
\end{theorem}

We can immediately conclude that our Theorem \ref{tt1} can be applied.

\begin{corollary}
Assume  (\ref{lip}), (\ref{cc4}), (\ref{u0}) and  (\ref{cc3})  and let $(u(t,x), t\in [0, T], x\in [0,1])$ be the unique solution to (\ref{burgers}). Then for every $t\in (0, T], x\in \mathbb{R}$, the random variable $u(t,x)$ admits a density  (denoted $f_{t, x}$) which satisfies
$$f_{t, x} \in \mathcal{B} _{1, \infty}^ {p} \mbox{ with } p< \alpha \wedge \frac{1}{2}.$$ 
\end{corollary}

Notice that the above result also applies to the stochastic heat equation with space variable in $[0,1]$, by choosing $g\equiv0$.

\subsection{The heat equation on $\mathbb{R}$}

Take $g\equiv  0$ and assume $I=\mathbb{R}$. The existence, the H\"older regularity and the moment estimates for the solution to (\ref{burgers}) have been obtained in \cite{Walsh}.  Actually, the relations (\ref{10f-9}) and (\ref{10f-10}) hold true in the case of the heat equation with spatial variable on the whole real line. The absolute continuity of the law of the solution is also well-known (see e.g. \cite{BP}).

In this case, our result in Theorem \ref{tt1} gives  in addition the Besov regularity of the solution. 

\begin{corollary}
Let $g\equiv 0$ and $I=\mathbb{R}$ in (\ref{burgers}) Assume (\ref{cc1}), (\ref{lip}),  (\ref{u0}), (\ref{cc3}) and (\ref{cc4}),  and let $(u(t,x), t\in [0, T], x\in \mathbb{R})$ be the unique solution to (\ref{burgers}). Then for every $t\in (0, T], x\in \mathbb{R}$, the random variable $u(t,x)$ admits a density  (denoted $f_{t, x}$) which satisfies
$$f_{t, x} \in \mathcal{B} _{1, \infty}^ {p} \mbox{ with } p< p< (\alpha \wedge \frac{1}{2}).$$ 
\end{corollary} 

We retrieve a result from \cite{Sanz1}, see Remark 3.4.

\end{document}